\documentclass[11pt,a4paper]{article}

\usepackage{bm,amsmath,amssymb, mathtools}
\usepackage{graphicx}
\usepackage{lineno,hyperref}
\usepackage{fancyhdr}
\usepackage[ruled]{algorithm2e}

\newcommand{\TFC}{Theory of Functional Connections}
\newcommand{\andd}{\qquad \text{and} \qquad}
\newcommand{\T}{^{\mbox{\tiny T}}}
\newcommand{\B}[1]{{\bm #1}}
\newcommand{\ds}{\displaystyle}
\DeclareMathOperator{\atantwo}{atan2}

\fancypagestyle{firststyle}
{
	\fancyhf{}
	\fancyfoot[LE,RO]{Page \thepage}
}

\begin{document}

\title{Bijective Mapping Analysis to Extend \\ the Theory of Functional Connections \\ to Non-rectangular 2-dimensional Domains}

\author{Daniele Mortari\thanks{Aerospace Engineering, Texas A\&M University, College Station TX, USA}, David Arnas\thanks{Massachusetts Institute of Technology (MIT), Cambridge, MA, USA.}}

\date{}	

\maketitle{} 	

\thispagestyle{firststyle}

\begin{abstract}
    This work presents an initial analysis of using bijective mappings to extend the \TFC\ to non-rectangular two-dimensional domains. Specifically, this manuscript proposes three different mappings techniques: a) complex mapping, b) projection mapping, and c) polynomial mapping. In that respect, an accurate least-squares approximated inverse mapping is also developed for those mappings having no closed-form inverse. Advantages and disadvantages of using these mappings are highlighted and a few examples are provided. Additionally, the paper shows how to replace boundary constraints expressed in terms of a piecewise sequence of functions with a single function, that is compatible and required by the \TFC\ already developed by rectangular domains.
\end{abstract}

\section{Introduction}

The Theory of Functional Connections (TFC) is a mathematical methodology to perform functional interpolation, that is, the process of deriving functionals, called \textit{constrained expressions}, which contain the constraints of the problem already embedded on their expression. To give an example, consider the following functional $y (x, g(x))$,
\begin{equation}\label{ce}
\begin{array}{rl}
y (x, g(x)) = g (x) &+ \quad \dfrac{2 g (-3) - 2 g (\pi) + (9 - \pi^2) (1 - g'(1))}{2\pi - \pi^2 + 15} \, x \\[10pt] ~ & + \quad \dfrac{g (\pi) - g (-3) + (\pi + 3) (1 - g'(1))}{2\pi - \pi^2 + 15} \, x^2,\end{array}
\end{equation}
where $x$ is the independent variable, $g$ is a function of $x$, and ($'$) indicates first derivative. This functional \textit{always} simultaneously satisfies the following two constraints,
\begin{equation*}
y (-3) = y(\pi) \andd y' (1) = 1,
\end{equation*}
{\it no matter what the $g (x)$ function is}! Function $g (x)$, called {\it free function}, can be any kind of function, including discontinuous or the Dirac functions, as long as {\it it is defined where the constraints are specified}. In this example the free function, $g (x)$, must be defined $x = 3$, at $x = \pi$, and it's derivative at $x =-3$. A mathematically proof has shown that constrained expressions, like the one given in Eq. (\ref{ce}), can be used to represent the whole set of functions satisfying the set of constraints they are derived for.

The TFC has been mainly developed \cite{TFC01,TFC09,TFC13,PDE01} to better solve constraint optimization problems, such as ODEs \cite{TFC02,TFC03,TFC04,TFC14}, PDEs \cite{PDE01, PDE02}, or programming \cite{TFC07,TFC08}, with effective applications in optimal control \cite{TFC05,TFC06,TFC16}, as well as in machine learning \cite{PDE01,TFC10,TFC11,TFC12}. In fact, a constrained expression restricts the whole space of functions to the subspace fully satisfying the constraints. This way, a constraint optimization problem can be solved using simpler, faster, more robust and accurate methods already developed and optimized to solve unconstrained optimization problems.

In this work we focus on developping mapping tools to extend the Theory of Functional Connections to non-rectangular domains in two dimensional spaces. This is done via domain mapping. In particular, three distinct domain mappings are analyzed. These are: a) {\it complex mapping}, b) {\it projection mapping}, and c) {\it polynomial mapping}. These three mappings are the subjects of the first three sections of this manuscript. These transformations are studied considering the direct and inverse mappings defined as:
\begin{enumerate}
	\item {\bf Direct mapping:} from rectangular domain $Z$ (coordinates $[a, b]$ if real, or $z = a + i \, b$, if complex with $a, b \in [-1, +1]$) to generic domain $W$ (coordinates $[x, y]$ if real, or $w = x + i \, y$, if complex).
	\item {\bf Inverse mapping:} from domain $W$ to domain $Z$. 
\end{enumerate}

The reason why this paper is specifically interested in the mapping from/to the $Z$ rectangular domain relies on the fact that the TFC has been fully developed for rectangular domains. Therefore, the main purpose of this study consists of developing easy mappings between the unit-square domain and any other domain, that flexibly can accommodate typical domains appearing in physics, science, and engineering problems. In other words, by obtaining these mappings, the TFC framework is consequently extended to generic domains.

In addition, this work extends the application of the TFC to discoutinuous boundary constraints, which appear in the projection mapping and may also appear in some other applications. This has been done by replacing the discontinuous sequence of boundary contraint functions with a single function. Finally, a least-squares method to approximate inverse transformations is provided for those mapping transformations with no analytical inverse.

\section{Complex (conformal) mapping}

In mathematics, a conformal mapping is a transformation that, locally, preserves angles, but not necessarily lengths\footnote{The conformal property is mathematically described in terms of the Jacobian of a coordinate transformation. If the Jacobian at each point is a positive scalar times a rotation matrix, then the transformation is defined as conformal.}. Since any function of a complex variable, $w = f (z)$, provides a complex variable, then a complex function can be seen as a complex change of variables or, from a different point of view, as a mapping in a complex Euclidean plane \cite{TFC15}. In particular, any complex mapping, $w = f (z)$, is mathematically proven to be conformal, that is, it is a trasformation preserving (signed) angles, meaning, local orientations. 

Conformal (angle-preserving) and area/volume-preserving mappings are particularly important transformations since they are very interesting for their applications in physics, science, and engineering. For instance, these mappings are used to transform differential equations in order to simplify the system and make it easier to solve \cite{TFC15}. 

A generic complex analytical mapping can be expressed as,
\begin{equation}
w = w_r + i \, w_i = f (z) = f (z_r + i \, z_i) = f_r (z_r, z_i) + i \, f_i (z_r, z_i),
\end{equation}
where $z,w \in \mathbb{C}$, $z_r, z_i, w_r, w_i \in \mathbb{R}$ represent the original and transformed points, and $f: \mathbb{C} \to \mathbb{C}$ is the transforming function, where  $f_r, f_i: \mathbb{R} \to \mathbb{R}$ are the real and imaginary parts of the transformation. Examples of some basic complex mappings are: 1) Translation, $w = z + z_t$, 2) Scaling, $w = \rho \, z$, 3) Rotation, $w = e^{i \, \phi} \, z$, 4) Affine, $w = \alpha \, z + \beta$, 5) Inversion, $w = 1/z$, 5) Exponential, $w = e^z$, 6), Squaring, $w = z^2$, 6) Cayley, $w = (z - 1)/(z + 1)$, and 7) M\"{o}bius (linear fractional), $w = (\alpha \, z + \beta)/(\gamma \, z + \delta)$, which is subject to $\alpha \, \delta \ne \beta \, \gamma$. In particular, the set of M\"{o}bius mapping includes, as a subset, translations, scalings, affines, inversions, and Cayley mappings. 
\begin{figure}[ht]
	\centering\includegraphics[width=\linewidth]{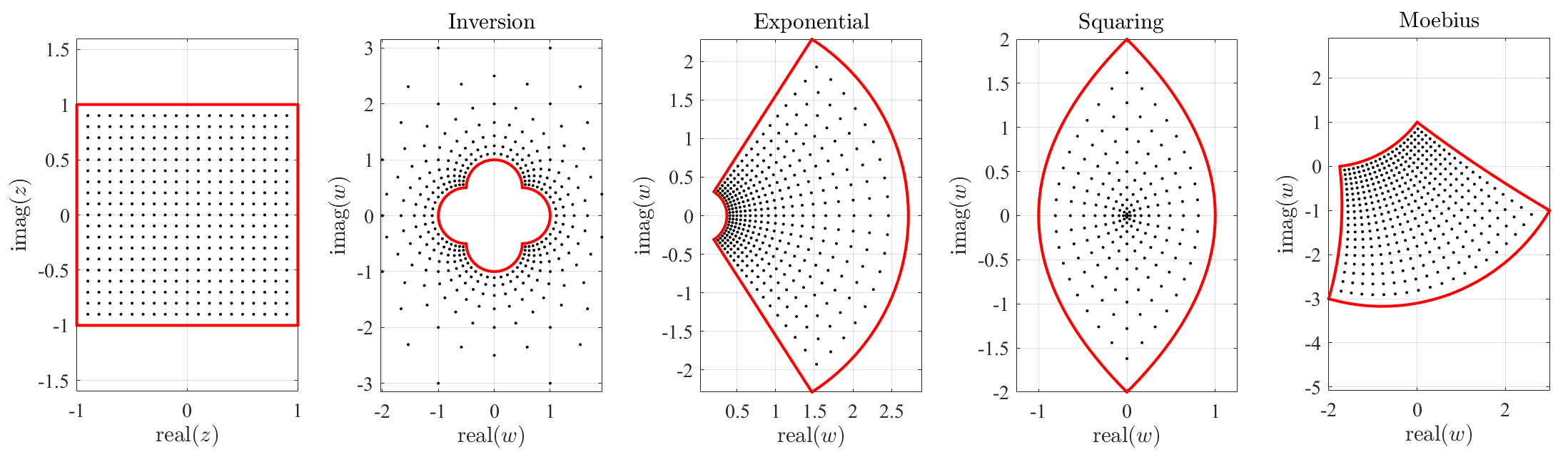}
	\caption{Examples of complex conformal mappings}
	\label{examples}
\end{figure}

Figure \ref{examples} shows the effect of inversion, exponenial, squaring, and M\"{o}bius mappings for $z = a + i \, b$ within the $a, b \in [-1, +1]$ domain. In particular, the M\"{o}bius mapping was selected with $\alpha =-3.8 + 4.4 \, i$, $\beta = 3.4 - 0.2 \, i$, $\gamma = 1$, and $\delta =-0.6 + 3.8 \, i$. These values were obtained by selecting $\gamma = 1$ and imposing the mapping
\begin{eqnarray*}
	z_1 =-1 - i & \qquad \to \qquad & w_1 = 2 - 3 \, i \\ 
	z_2 = +1 - i& \qquad \to \qquad & w_2 = 3 - i \\ 
	z_3 = +1 + i& \qquad \to \qquad & w_3 = i
\end{eqnarray*}
allowing the computation of $\alpha$, $\beta$, and $\delta$. The M\"{o}bius mapping can be seen as the stereographic projection from the plane to a rotated and translated unit-sphere. The M\"{o}bius mapping plays a particularly important role in the conformal complex mapping because: 1) it admits a closed-form expression for the inverse mapping and 2) is the composition or the most basic transformations of the space, including similarities, orthogonal transformations, and inversions. Also and for a number of dimensions equal or higher than 3, Liouville's theorem states that M\"{o}bius transformations are the unique transformations that are conformal.

\subsection{Control points in complex conformal mapping}

\begin{figure}[!ht]
	\centering\includegraphics[width=\linewidth]{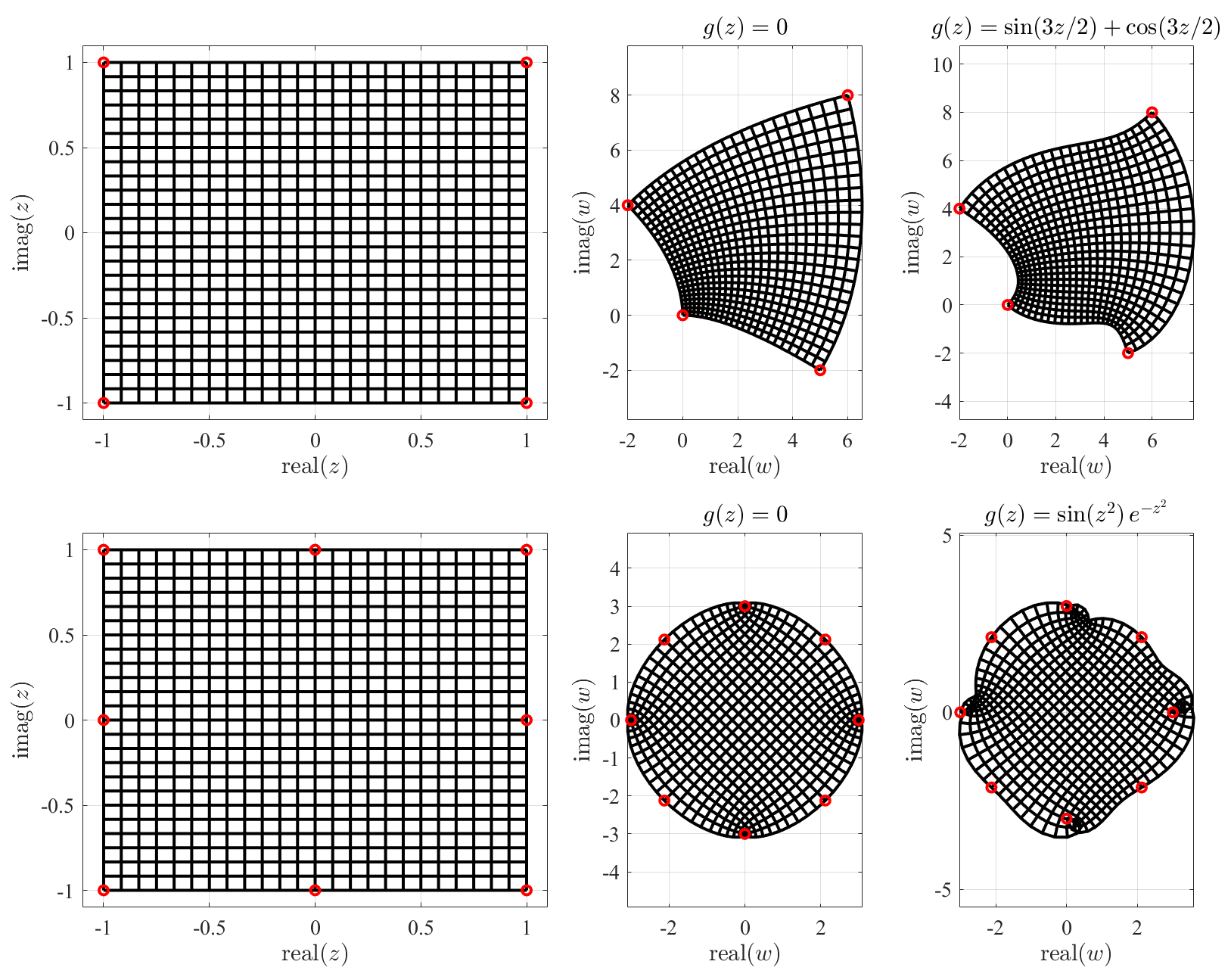}
	\caption{Rubber effect on the $W$ domain by the TFC free function.}
	\label{Fig3}
\end{figure}

Given a set of $n$ complex control points, $w_i$ with $i \in\{1,\dots,n\}$, in the $W$ domain associated to $n$ corresponding boundary reference points, $z_j$ with $j \in\{1,\dots,n\}$, in the unit-square $Z$ domain, the transformation,
\begin{equation}\label{eq02}
w (z) = w_i \, \phi_i (z) \qquad \text{where} \qquad \phi_i (z_j) = \delta_{ij}
\end{equation}
is a $n$-point complex conformal mapping. In this equation, $\phi_i (z)$ are a set of switching functions related to the reference points, $z_j$. In particular, Let select $n = 4$ control points (red markers on top plots of Fig. \ref{Fig3}) defined by,
\begin{equation}\label{points}
\text{reference points:} \; \begin{cases} z_1 =-1 - i\\ z_2 = +1 - i\\ z_3 = +1 + i\\ z_4 =-1 + i\end{cases} \quad \text{control points:} \; \begin{cases} w_1 = 0 \\ w_2 = +5 - 2 \, i\\ w_3 = +6 + 8 \, i\\ w_4 =-2 + 4 \, i\end{cases}.
\end{equation}
The switching functions, $\phi_i (z)$, are expressed as a linear combination of $n$ complex support functions. For the sake of simplicity, let the support functions be  selected as complex monomials. The main reason for that is the simple closed-form expressions of derivatives and integrals for this set of functions. Note that, since {\it any} set of $n$ linearly independent functions can be selected as support functions, this provides additional degrees of freedom that are not analyzed in this study. Therefore, using monomial, the switching functions are expressed as,
\begin{equation*}
\phi_i (z) = c_{ij} \, z^{j-1}
\end{equation*}
where $c_{ij}$ are a set of constants. Since the switching property states that $\phi_i (z_j) = \delta_{ij}$, we can derive the following expression,
\begin{equation*}
\begin{bmatrix} c_{11} & c_{12} & c_{13} & c_{14} \\ c_{21} & c_{22} & c_{23} & c_{24} \\ c_{31} & c_{32} & c_{33} & c_{34} \\ c_{41} & c_{42} & c_{43} & c_{44}\end{bmatrix} = \begin{bmatrix} 1 & 1 & 1 & 1\\ z_1 & z_2 & z_3 & z_4\\ z_1^2 & z_2^2 & z_3^2 & z_4^2\\ z_1^3 & z_2^3 & z_3^3 & z_4^3\end{bmatrix}^{-1} = \dfrac{1}{16} \begin{bmatrix*}[r] 4 & -2 + 2 i & -2 i & 1 - i \\ 4 & 2 + 2 i & 2 i & -1 + i \\ 4 & 2 - 2 i & -2 i & -1 - i \\ 4 & -2 + 2 i & 2 i & 1 - i\end{bmatrix*},
\end{equation*}
which provides the coefficients needed for the direct (from $Z$ to $W$) mapping given in Eq. (\ref{eq02}). 

It is important to outline that in general complex mapping may not be bjiective, which makes impossible to find an inverse transformation for the whole domain. However, due to the switching functions, the inverse mapping always exists for the control points. When the complex mapping is bijective, then a least-squares approximated inverse mapping using a basis of orthogonal polynomials and presented in Section~\ref{ApproxInverse}, can be used.

\subsection{Complex TFC mapping}

We can generalize the previous result by using the \TFC. Particularly, the functional,
\begin{equation}\label{eq02b}
w (z, g (z)) = g (z) + \ds\sum_{i = 1}^n 
[w_i - g (z_i)] \, \phi_i (z),
\end{equation}
where $g (z)$ is a free complex function, can be seen as the most general expression that defines the conformal mapping between the control points in $Z$ and $W$. In fact, Eq.~\eqref{eq02b} represents the TFC functional generalization of Eq. (\ref{eq02}).

Figure~\ref{Fig3} shows one effect of the free function in the complex mapping defined by the control points from Eq.~\eqref{points}. The free function modifies the domain while the conformal property is still preserved. The top-left plots of Fig. \ref{Fig3} shows a grid of orthogonal lines in the $Z$ domain, while the top-center plot shows the lines mapped in the corresponding $W$ domain using $g (z) = 0$, that is, using Eq. (\ref{eq02}), and the control points given in Eq. (\ref{points}). The free function does not change the mapping of the control points. These points are constraints for the $W$ domain. The top-right plot shows the mapped $W$ domain obtained using the free function $g (z) = \sin(3 z/2) + \cos(3 z/2)$, that is, applying Eq. (\ref{eq02b}). This shows the free function ``rubber'' effect on the $W$ domain. The bottom plots are obtained for 8 control points distributed along a circle and using for the free function, the expression $g (z) = \sin(z^2) \, e^{-z^2}$. This example highlights that the selection of control points is not completely free because they often produce non-bjiective mapping (points in the $W$ domain may be associated to multiple points in the $Z$ domain).

\section{Projection mapping}

In this section a transformation method based on a projection that preserves the continuity of the space and provides both, the direct and the inverse mapping, is presented. However, and due to its nature, the conformal property is not maintained by this mapping approach.

Projection mapping is based on the selection of a set of points (the projection points) inside the boundaries of the domain in such a way that a straight segment can be defined between any point in the domain and at least one projection point, where this segment does not cross at any moment any of the boundaries of the domain. This means that the distances from these projection points to the points at the boundary can be expressed as a continuous function. The most interesting feature of the projection mapping is that, using an identical formalism, the direct and the inverse mapping transformation can be derived (isomorphism). If the projection points are selected properly, the continuity of the transformation can be assured, however, in general, the derivative of this transformation is not defined in the whole domain and thus, the transformation is not conformal.

In order to show the methodology in a clearer manner, we first present the case of a single projection point. Then, we use this result to generalize the methodology for the case of muliple projection points. 

\subsection{Single-point projection mapping}

The idea of this methodology is to perform a linear transformation in polar coordinates with respect to a projection point that lies inside the domain. While the projection mapping can be applied to any dimensional domains, to be consistent with the complex conformal mapping, the analysis is here restricted to bivariate domains. 

Let $\{r, \theta\}$ be the polar coordinates in the space $Z$, and let $\{d, \phi\}$ the polar coordinates in the space $W$. We know that for any point inside the square boundary of $Z$, the boundary can be defined as $r_b(\theta) = f(\theta)$ where $f$ is a continuous but not differentiable function. In the same way, if the boundaries in $W$ can be expressed as $d_b(\phi) = g(\phi)$ from a given projection point, then, the following direct transformation can be defined:
\begin{eqnarray}
\phi & = & \theta; \nonumber \\
d & = & d_b(\theta)\displaystyle\frac{r}{r_b(\theta)};
\end{eqnarray}
whose inverse transformation is:
\begin{eqnarray*}
	\theta & = & \phi; \nonumber \\
	r & = & r_b(\phi)\displaystyle\frac{d}{d_b(\phi)}.
\end{eqnarray*}

Nevertheless, in some applications it is useful to define the boundary in $W$ as a given polygon defined by a set of control points. Therefore, we include in this section the algorithm to generate a transformation mapping between a set of $n$ control points. In this case, and since we are using the control points as the reference to perform the transformation, the angles in both mappings are not equal as in the previous equations. Particularly, a linear interpolation in the angle between the control points and the projection point is performed in order to maintain the continuity of the transformation.
\label{alg}
\begin{algorithm}[!ht]
	\caption{Direct transformation by projection mapping using a set of $n$ control points}
	Angle and distance with respect to the projecting point in $Z$: \\
	$\theta = \atantwo(b - c_b, a - c_a);$
	$r^2 = (a - c_a)^2 + (b - c_b)^2;$ \\
	\For{$j = 1$ to $j = n$}{
		Angle boundaries of the region in $Z$: \\
		$\theta_{min} = \atantwo(p_b (j) - c_b, p_a (j) - c_a);$ 
		$\theta_{max} = \atantwo(p_b (j + 1) - c_b, p_a (j + 1) - c_a);$
		$\Delta\theta = \theta_{max}-\theta_{min} \mod(2\pi);$ \\
		\If {$((\theta \geq \theta_{min}) \text{ and } (\theta \leq \theta_{max})) \text{ or } ((\theta \geq \theta_{min} \mod(2\pi)) \text{ and } (\theta \leq \theta_{max} \mod(2\pi)))$}{
			Line passing through the projecting point in $Z$: \\
			$m_p = (b - c_b)/(a - c_a);$
			$n_p = c_b - m_pc_a;$ \\
			Line defining the boundary in $Z$: \\
			$m_b = (p_b(j+1) - p_b(j))/(p_a(j+1) - p_a(j));$
			$n_b = p_b(j) - m_b*p_a(j);$ \\
			Distance to the intersecting point in $Z$: \\
			$r_b^2 = ((n_p-n_b)/(m_b-m_p)-c_a)^2 + (m_b(n_p-n_b)/(m_b-m_p)+n_b-c_b)^2;$ \\
			Angle boundaries of the region in $W$: \\
			$\phi_{min} = \atantwo(q_y(j)-\upsilon_y,q_x(j)-\upsilon_x);$
			$\phi_{max} = \atantwo(q_y(j+1)-\upsilon_y,q_x(j+1)-\upsilon_x);$
			$\Delta\phi = \phi_{max}-\phi_{min} \mod(2\pi);$ \\
			Transformed angle: \\
			$\phi = \theta-\theta_{min} \mod(2\pi);$
			$\phi = \phi\Delta\phi/\Delta\theta + \phi_{min};$ \\
			Line passing through the projecting point in $W$: \\
			$m_p = \tan(\phi);$
			$n_p = \upsilon_y - \upsilon_xm_p;$ \\
			Line defining the boundary in $W$: \\
			$m_b = (q_y(j+1) - q_y(j))/(q_x(j+1) - q_x(j));$
			$n_b = q_y(j) - m_bq_x(j);$ \\
			Distance to the intersecting point in $W$: \\
			$d_b^2  = ((n_p-n_b)/(m_b-m_p)-\upsilon_x)^2 + (m_b(n_p-n_b)/(m_b-m_p)+n_b-\upsilon_y)^2;$ \\
			break; \\
		}
	}
	Transformated point: \\
	$d = d_br/r_b;$ \\
	$x = \upsilon_x + d\cos(\phi);$
	$y = \upsilon_y + d\sin(\phi);$
\end{algorithm}

We define $a$ and $b$ to the coordinates in $Z$ of a given point inside the original domain that is required to be transformed, and $x$ and $y$ to the coordinates in $W$ of the transformed point. Moreover, $\B{c} = \{c_a, c_b\}$, $\B{\upsilon} = \{\upsilon_x, \upsilon_y\}$ are the coordinates of the projections points in $Z$ and $W$ respectively, while $\{p_a (j), p_b (j)\}$, and $\{q_x (j), q_y (j)\}$, with $j\in\{1, \dots, n + 1\}$, are the coordinates of the $n$ control points defined in $Z$ and $W$ respectively plus the first one repeated in order to close the polygon. Note that the control points have to be consequtive, that is, it must be possible to generate the polygon by linking the control points in the order provided. Then the transformation of $\{a, b\}$ into $\{x, y\}$ can be performed using Algorithm~\ref{alg}. Finally, the inverse transformation is equivalent to the algorithm presented.
\begin{figure}[!ht]
	\centering\includegraphics[width=\linewidth]{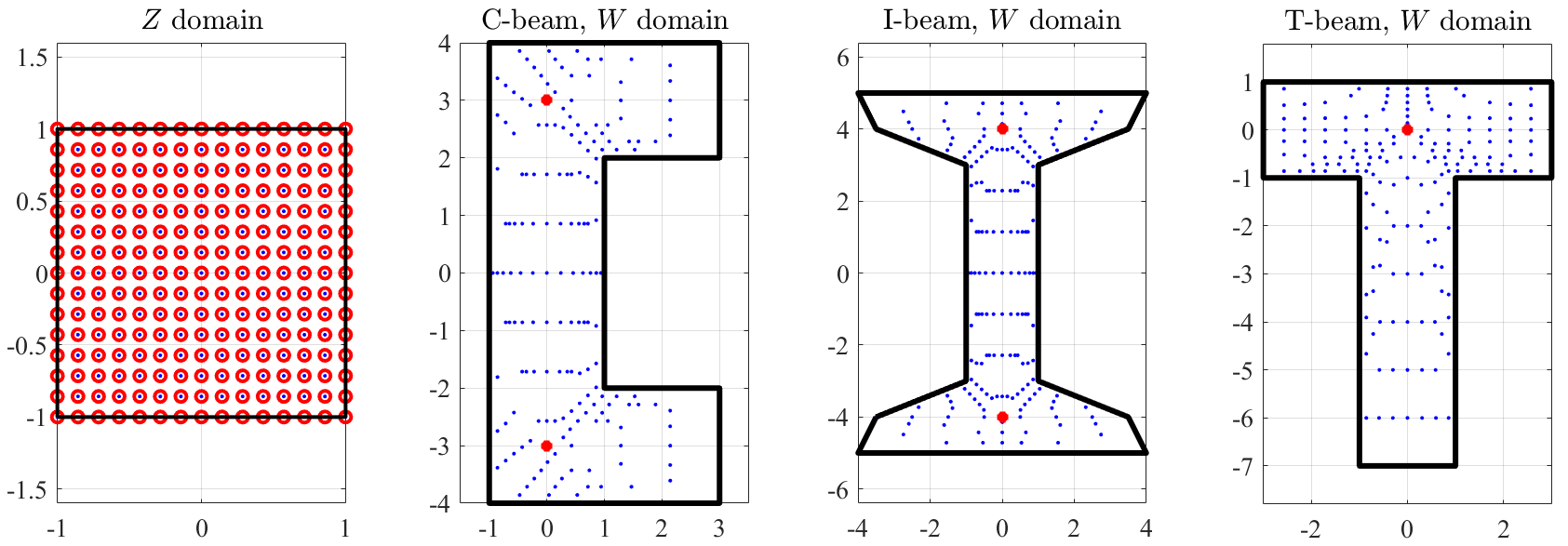}
	\caption{Examples of projection mapping for ``C,'' ``I,'' and ``T'' beam shape domains}
	\label{beams}
\end{figure}

As an example of this method, Fig.\ref{beams} (right) presents the result of applying the transformation (direct and inverse) to the control points provided by Table~\ref{tab:tcp} and the projecting points $c = \{0, 0\}$ and $\upsilon = \{0, 0\}$. As it can be seen, the polynomial in $W$ resembles a ``T'' shape, where each consecutive pair of control points generates a region where the transformation is continuous and differentiable. However, at the boundaries between regions the transformation stops to be differentiable, although it mantains the property of continuity.

\begin{table}[!ht]
	\centering
	\begin{tabular}{|c|c|c|c|c|c|c|c|c|c|}
		\hline
		$j$ & 1 & 2 & 3 & 4 & 5 & 6 & 7 & 8 & 9 \\
		\hline
		$p_{a}$ & -1 & -0.5 & 0 & 0.5 & 1 & 1 & 1 & 1 & 1 \\
		$p_{b}$ & -1 & -1 & -1 & -1 & -1 & -0.5 & 0 & 0.5 & 1 \\
		\hline
		$q_{x}$ & -2 & -1 & -1 & -1 & 0 & 1 & 1 & 1 & 2 \\ 
		$q_{y}$ & -1 & -1 & -4 & -7 & -7 & -7 & -4 & -1 & -1 \\
		\hline
		\hline
		$j$ & 10 & 11 & 12 & 13 & 14 & 15 & 16 & 17 & \\
		\hline
		$p_{a}$ & 0.5 & 0 & -0.5 & -1 & -1 & -1 & -1 & -1 & \\
		$p_{b}$ & 1 & 1 & 1 & 1 & 0.5 & 0 & -0.5 & -1 & \\
		\hline
		$q_{x}$ & 3 & 3 & 3 & 0 & -3 & -3 & -3 & -2 & \\ 
		$q_{y}$ & -1 & 0 & 1 & 1 & 1 & 0 & -1 & -1 & \\
		\hline
	\end{tabular}
	\caption{Control points for the ``T'' beam shape transformation}
	\label{tab:tcp}
\end{table}

\subsection{Multiple-point projection mapping}

There are cases in which it is not possible to generate a function $d_b(\phi) = g(\phi)$ as in the previous case due to the shape of the domain in $W$. In these situations, it is possible to solve these problems using a set of projection points instead of just one, however, some considerations regarding the continuity of the transformation must to be considered first. In that regard, and as a first approach to deal with this situation, we could think on finding the smallest number of projection points required to cover the whole domain. This represents in effect the so called ``art gallery problem''.

The ``art gallery problem'' (known also as the museum problem) is a problem in computational geometry providing the solution for the minimum number of required points to cover the complete domain. Its name comes from the real-world problem of placing the minimum number of guards in an art gallery such that, all together, can observe the whole gallery. This problem was first solved by Chv\'{a}tal in 1975 \cite{Chvatal}, but Fisk, in 1978, provided such a brilliant and short proof \cite{Fisk} that Paul Erd\H{o}s decided to include it in ``The BOOK,'' an imaginary book in which God keeps the most elegant proof of each mathematical theorem. This book was then published after Erd\H{o}s dead \cite{BOOK}. The solution to this problem states that $\left\lfloor n/3 \right\rfloor$ guards are always sufficient and (sometimes) necessary to guard a simple polygon with $n$ vertices. In general these points are in the boundary of the polygon and thus, although the theorem is interesting, it is not really practical for most applications since it distributes points more densily in the boundaries of the domain in $W$ rather than obtaining a more uniform distribution of points.

Therefore, instead of focusing on identifying the smallest number of projection points required, in this work we deal with the problem on how to separate both domains (in $Z$ and $W$) in a set of regions, one for each projecting point, in such a way that the property of continuity is preserved through the transformation. Once these regions are defined, we proceed as in the case of single-point projection mapping. Note that by doing that, each time that a point is transformed, the method requires to determine the region, from the ones defined, to which it belongs. This can be done, for instance, by using the transformation for each of these regions chequing if the point is inside the region in the $Z$ space.

\begin{figure}[!ht]
	\centering\includegraphics[width=0.7\linewidth]{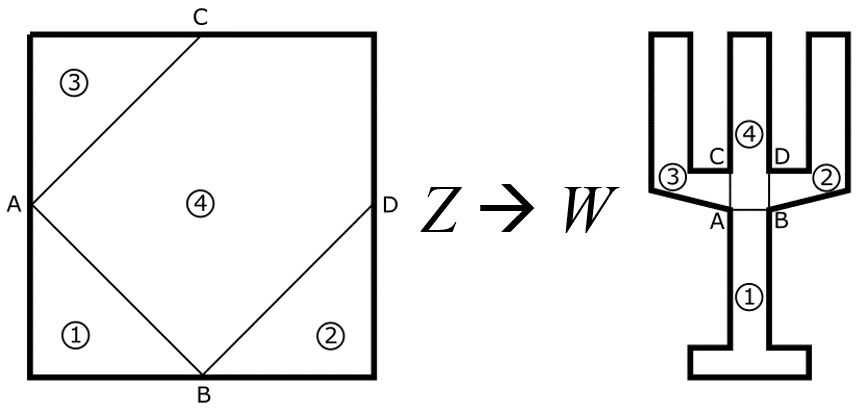}
	\caption{Definition on region for multi-point projection mapping}
	\label{multi}
\end{figure}

Let $\delta\Omega$ be the boundary in $W$ of a given domain. Then, connecting with a line two points in $\delta\Omega$ implies connecting two points in the boundary of $Z$ since during the transformation the boundary in one space transforms into the boundary on the other space. For instance, from Figure~\ref{multi}, if we select the points A and B from $W$, we know that both points belong to the boundary in $Z$ and also, both points can be connected through a curve in $Z$ and $W$. Then, if we use just straight lines as curves, this means that A and B in $Z$ must be in a different side of the square boundary. The same can be said with the points $C$ and $D$. Following this procedure, it is possible to generate 4 regions in which each projection point can cover completely its own region. At this point, it is important to note that the boundary between regions 1 and 4 is the curve between A and B in both spaces. In the same way, regions 3 and 4 have a boundary in the curve between A and C, while in regions 2 and 4 the boundary in common is defined by the curve between B and D. This means that since the transformation is continuous inside the regions and at the boundaries of these regions, then, the whole transformation is continuous in the whole domain. This process can be continued even using the boundaries between these defined regions. This means that, for instance, if region 4 is split in half vertically, a point in the curve between A and B will be part of the new boundary created between regions 4' and 4'', while another point in the boundary should be in $\delta\Omega$ between C and D. 

In order to present an example of application for the case of multiple projection points, we apply this methodology to the ``C'' and ``I'' beam transformations shown in Fig.~\ref{beams} . To that end, we know that for these cases we can cover all the domain with two projection points. Moreover, and since both figures are symmetric with respect $b = 0$ and $y = 0$ repectively, we use this line as the boundary betwwen the two regions in which to apply Algorithm~\ref{alg}. In particular, for the case of the ``C'' beam shape, we define the control points for the two defined regions in Table~\ref{tab:ccp}, where the projecting points used are $\B{c_1} = \{0, 0.5\}$, $\B{c_1} = \{0, -0.5\}$, $\B{\upsilon_1} = \{0, 3i\}$, and $\B{\upsilon_2} = \{0, -3i\}$. Figure~\ref{beams} shows the result of the transformation.

\begin{table}[!ht]
	\centering
	\begin{tabular}{|c|c|c|c|c|c|c|c|c|c|c|}
		\hline
		$j$ & 1 & 2 & 3 & 4 & 5 & 6 & 7 & 8 & 9 & 10 \\ \hline\hline
		$p_a'$ & 1 & 1 & 1 & 0.5 & 0 & $-0.5$ & $-1$ & $-1$ & $-1$ & 1 \\
		$p_b'$ & 0 & 0.5 & 1 & 1 & 1 & 1 & 1 & 0.5 & 0 & 0 \\[2pt]
		\hline
		$p_a''$ & $-1$ & $-1$ & $-1$ & $-0.5$ & 0 & 0.5 & 1 & 1 & 1 & $-1$ \\
		$p_{b}''$ & 0 & $-0.5$ & $-1$ & $-1$ & $-1$ & $-1$ & $-1$ & $-0.5$ & 0 & 0 \\[2pt]
		\hline\hline
		$q_x'$ & 1 & 1 & 2 & 3 & 3 & 3 & 1 & $-1$ & $-1$ & 1 \\
		$q_y'$ & 0 & 2 & 2 & 2 & 3 & 4 & 4 & 4 & 0 & 0 \\[2pt]
		\hline
		$q_x''$ & $-1$ & $-1$ & 1 & 3 & 3 & 3 & 2 & 1 & 1 & $-1$ \\
		$q_y''$ & 0 & $-4$ & $-4$ & $-4$ & $-3$ & $-2$ & $-2$ & $-2$ & 0 & 0 \\[2pt]
		\hline
	\end{tabular}
	\caption{Control points for the ``C'' beam shape transformation}
	\label{tab:ccp}
\end{table}

\subsection{Maintaining the density of points during the transformation}

In some applications it is of interest to maintain the density of points through the transformation between $W$ and $Z$. This can be achieved easily using the projection mapping transformation proposed by defining appropriately  the control points in the square domain in $Z$. 

Let $n$ be the number of uniformly distributed points in $W$ that are inside the problem domain of area $A$. Let $n_p$ be the number of points in a given region in $W$ of area $A_w$ that is defined by two consecutive control points and the projection point. Then, since the points are uniformly distributed: $n_p/n = A_w/A$. This argument can be also applied to the square domain in $Z$ leading to: $n_p/n = A_z/A_s$ where $A_z$ is the area defined by the transformed control points and the projection point in $Z$, and $A_s$ is the area of the square. This means that, since $n_p/n$ is a fixed value in the transformation, we have to impose that $A_z = A_s \, A_w/A$ in order to maintain the density of points. In other words, the proportion of the area or each region with respect with the area of the domain must be equal in both domains. Therefore, it is poisible to select the two control points in the boundary of $Z$ in such a way that this area is maintained. For instance, if the area of a given region in $W$ is 10\% of the domain area, and one of the control points in $Z$ was already in the position $\{1,0\}$, this means that the second control point for that region should be $\{1,0.8\}$. It is important to note that in general, these new control points will not coincide with the corners of the square, and thus, additional control points in the square (and their transformed equivalents in $W$) must be generated in order to define completely the domain in $Z$.

\subsection{Merging a sequence of boundary functions into a single function}

The strength of the projection mapping consists of providing the direct and inverse transformations with the same formalism and with no approximation. However, for generic polygon domains this leads the TFC problem to be subject to boundary conditions made of a sequence of piecewise functions that in general can even be discontinuous. In this subsection using the Heaviside step function,
\begin{equation}
H (x) = \begin{cases} = 0 \qquad \text{if} \qquad x < 0 \\ = 1 \qquad \text{if} \qquad x > 0\end{cases}, \nonumber
\end{equation}
we show how to replace, with no approximation, a set of $n$ contiguous piecewise functions by a single function that can be used in the TFC framework. In particular, the possibility of discontinuities between subsequent piecewise functions (no $C^{\, 0}$ continuity) is also included.

\begin{figure}[!ht]
	\centering\includegraphics[width=\linewidth]{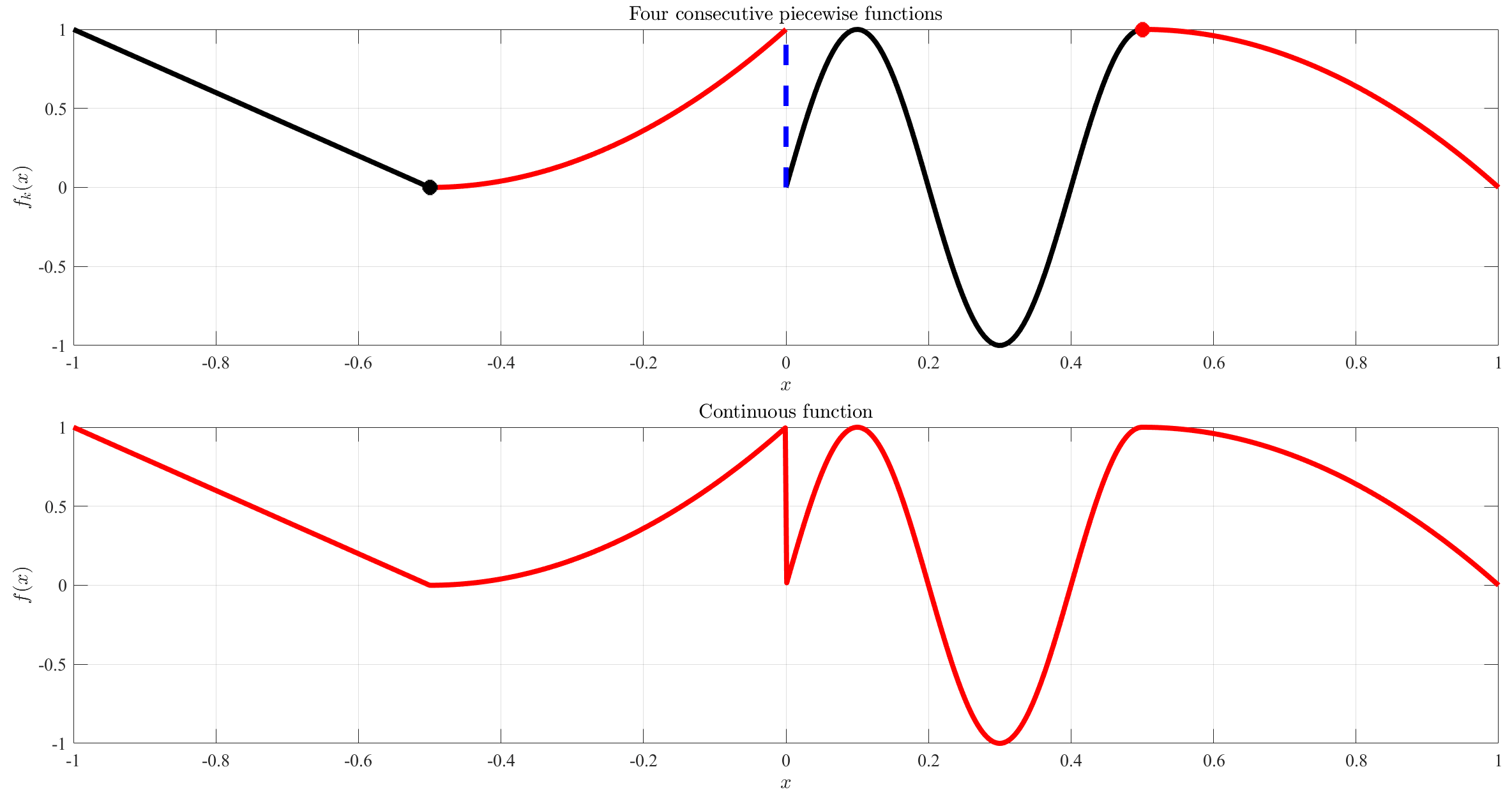}
	\caption{Four functions given in Eq. (\ref{pts}) (top) and equivalent single function given by Eq. (\ref{d02}) (bottom)}
	\label{fig:test1}
\end{figure}

Consider, for example, the boundary constraints described by the following set of four piecewise functions, each defined in contiguous ranges,
\begin{equation}\label{pts}
\begin{array}{lcl}
f_1 (x) =-1 - 2 \, x & \text{in} & x \in [-1, -1/2] \\
f_2 (x) = 1 + 4 (1 + x) \, x & \text{in} & x \in [-1/2, \; 0] \\
f_3 (x) = \sin(5\pi \, x) & \text{in} & x \in [0, \; +1/2] \\
f_4 (x) = 4 (1 - x) \, x & \text{in} & x \in [+1/2, +1]
\end{array}
\end{equation}
This piecewise function is shown in the top plot of Fig. \ref{fig:test1}; while the single function, $f (x)$, which is derived by the general equation,
\begin{equation}\label{d02}
\boxed{ f (x) = f_1 (x) + \ds\sum_{k = 2}^n H (x - x_k) \, \left[f_k(x) - f_{k-1} (x)\right] }
\end{equation}
applied for $n = 4$ piecewise functions, is shown in the bottom plot of Fig. \ref{fig:test1}. In particular, this example includes, $C^{\,0}$ continuity (black marker) at $x_2 =-1/2$, discontinuity (blue dashed line) at $x_3 = 0$, and $C^{\,0}$ and $C^{\,1}$ continuity (red marker) at $x_2 = 1/2$.
\begin{figure}[!ht]
	\centering\includegraphics[width=\linewidth]{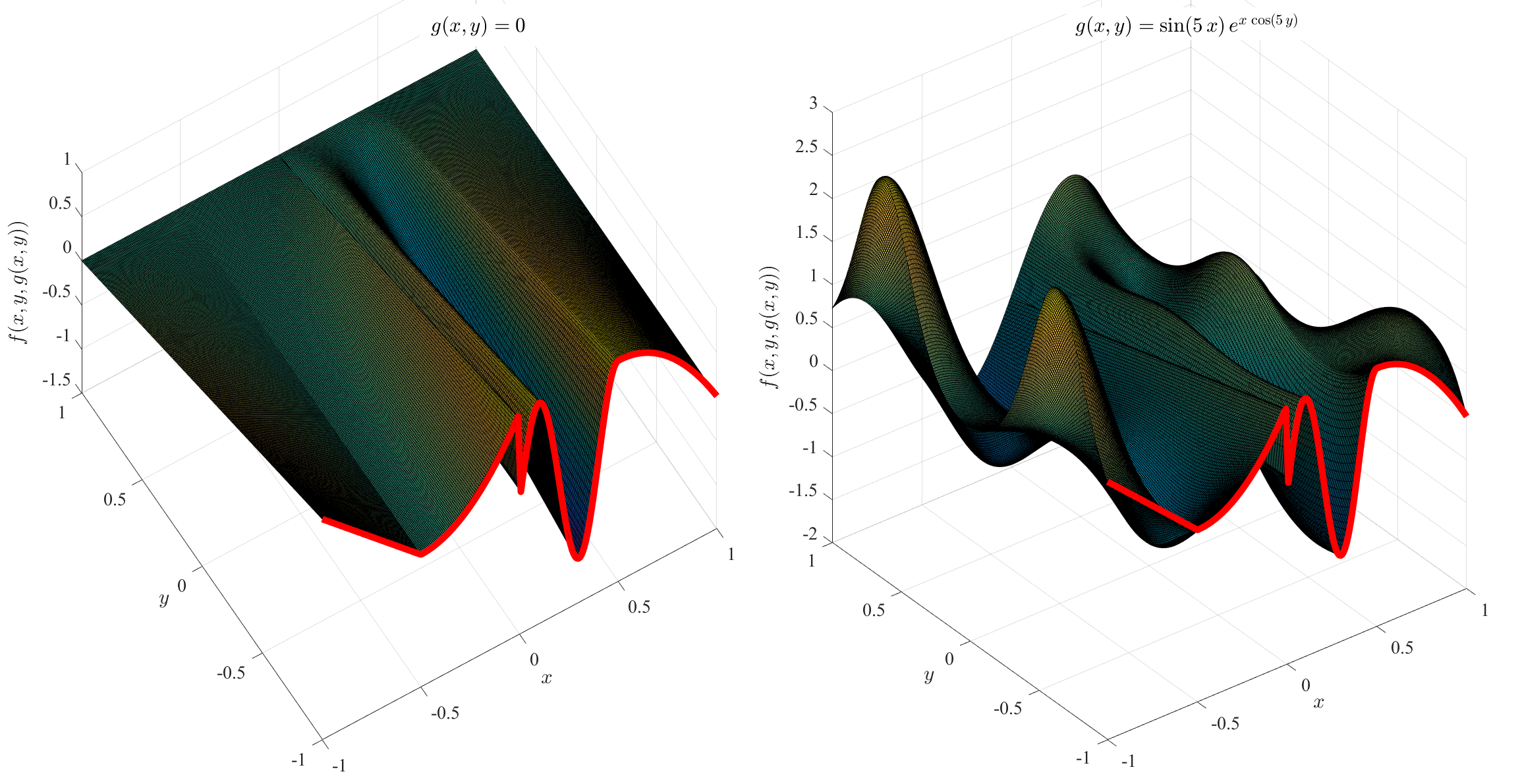}
	\caption{TFC surface examples obtained with $g (x) = 0$ (left) and $g (x) = \sin(5 \, x) \, e^{x \, \cos (5 \, x)}$ (right)}
	\label{fig:test2}
\end{figure}

The function, $f (x)$, computed using Eq. (\ref{d02}) can be used, for instance, to describe a boundary constraint at $y =-1$ in the 2-dimensional TFC matrix formulation \cite{TFC09,TFC13},
\begin{equation}
f (x, y, g (x, y)) = g (x, y) + \B{v}\T (x) \, \left[\mathcal{M} (f (x)) - \mathcal{M} (g (x, y))\right] \, \B{v} (y)
\end{equation}
where the matrix tensor is simply,
\begin{equation}
\mathcal{M} (f (x)) = \begin{bmatrix} 0, & f (x)\end{bmatrix}
\end{equation}
and the switching vectors are,
\begin{equation}
\B{v} (x) = 1 \andd \B{v} (y) = \dfrac{1}{2} \begin{Bmatrix} 2 \\ 1 - y\end{Bmatrix}
\end{equation}
Therefore, the TFC functional is,
\begin{equation}
f (x, y, g (x, y)) = g (x, y) + \dfrac{1 - y}{2} \, \left\{f (x) - g (x,-1)\right\}
\end{equation}
The single function, $f (x)$, can be used in the TFC framework as a boundary constraint. Figure \ref{fig:test2} shows two TFC surfaces using the single constraint boundary function, $f (x)$, shown in re. This has been done for two distinct expressions of the free functions: $g (x, y) = 0$ (left figure, most simple interpolating surface) and $g (x) = \sin(5 \, x) \, e^{x \, \cos (5 \, x)}$ (right figure).

\section{Polynomial mapping}

The polynomial mapping presented in this section is a mapping between two real domains. This mapping is performed using a set of $n$ polynomial (switching) functions and a set of $n$ boundary (control) points in the $Z$ and $W$ domains. While the proposed mapping can be applied to any dimensional domains, to be consistent with the complex conformal mapping, the analysis is then restricted to the bi-variate domains only. 

Let the boundaries of two domains, $Z$ and $W$, be identified by the following sequence of $n$ control points,
\begin{equation}
\B{z}_k = \begin{Bmatrix} a_k \\ b_k\end{Bmatrix} \andd \B{w}_k = \begin{Bmatrix} x_k \\ y_k\end{Bmatrix} \qquad \text{where} \quad k \in \{1,\dots, n\}
\end{equation}
then, the direct mapping is provided by,
\begin{equation}\label{switching}
\begin{Bmatrix} x \\ y\end{Bmatrix} = \ds\sum_{i = 1}^n \begin{Bmatrix} x_i \\ y_i\end{Bmatrix} \, f_i (a, b),
\end{equation}
where $f_i (a_j, b_j) = \delta_{ij}$ defines the switching property of the polynomials. This property gives the points association relationship, $[x_i, y_i] \leftrightarrow [a_i, b_i]$. The $n$ coefficients of the switching functions are computed from the switching function property. For instance, consider the simplest example of $n = 3$ points (triangle) mapping. In this case, the mapping becomes linear, $f_i (a, b) = c_{i1} + c_{i2} \, a + c_{i3} \, b$. The coefficients of the three switching functions ($i = 1, 2, 3$) are derived by imposing the switching functions property,
\begin{eqnarray}
& \begin{cases} 1 = c_{11} + c_{12} \, a_1 + c_{13} \, b_1 \\ 0 = c_{11} + c_{12} \, a_2 + c_{13} \, b_2 \\ 0 = c_{11} + c_{12} \, a_3 + c_{13} \, b_3\end{cases} \quad \begin{cases} 0 = c_{21} + c_{22} \, a_1 + c_{23} \, b_1 \\ 1 = c_{21} + c_{22} \, a_2 + c_{23} \, b_2 \\ 0 = c_{21} + c_{22} \, a_3 + c_{23} \, b_3\end{cases} & \nonumber \\ 
& \text{and} \; \begin{cases} 0 = c_{31} + c_{32} \, a_1 + c_{33} \, b_1 \\ 0 = c_{31} + c_{32} \, a_2 + c_{33} \, b_2 \\ 1 = c_{31} + c_{32} \, a_3 + c_{33} \, b_3\end{cases} &
\end{eqnarray}
or, in matrix form,
\begin{eqnarray}
& \begin{bmatrix} c_{11} & c_{12} & c_{13} \\ c_{21} & c_{22} & c_{23} \\ c_{31} & c_{32} & c_{33}\end{bmatrix} \begin{bmatrix} 1 & 1 & 1 \\ a_1 & a_2 & a_3 \\ b_1 & b_2 & b_3\end{bmatrix} = \begin{bmatrix} 1 & 0 & 0 \\ 0 & 1 & 0 \\ 0 & 0 & 1\end{bmatrix} \quad \to & \nonumber \\ 
& \to \quad \begin{bmatrix} c_{11} & c_{12} & c_{13} \\ c_{21} & c_{22} & c_{23} \\ c_{31} & c_{32} & c_{33}\end{bmatrix} =  \begin{bmatrix} 1 & 1 & 1 \\ a_1 & a_2 & a_3 \\ b_1 & b_2 & b_3\end{bmatrix}^{-1} &
\end{eqnarray}
\begin{figure}[ht]
	\centering\includegraphics[width=\linewidth]{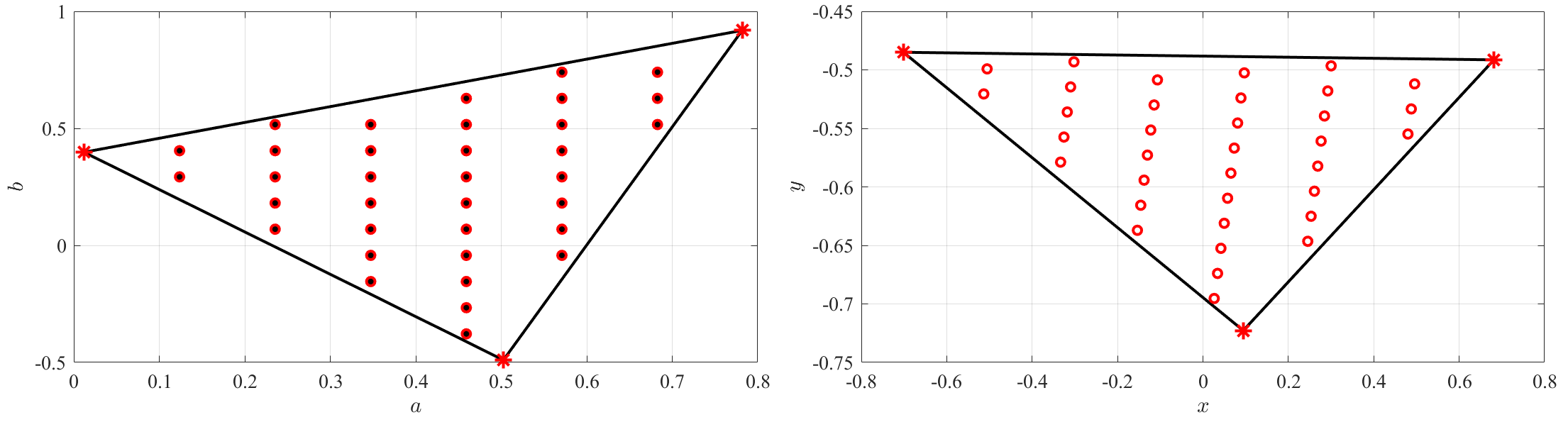}
	\caption{3-points polynomial mapping example}
	\label{triangle}
\end{figure}

In this example the mapping function is linear and, therefore, the inverse mapping function,
\begin{equation}
\begin{Bmatrix} a \\ b\end{Bmatrix} = \ds\sum_{i = 1}^n \begin{Bmatrix} a_i \\ b_i\end{Bmatrix} \, p_i (x, y)
\end{equation}
where
\begin{equation}
p_i (x, y) = c_{i1} + c_{i2} \, x + c_{i3} \ y \quad \to \quad \begin{bmatrix} c_{11} & c_{12} & c_{13} \\ c_{21} & c_{22} & c_{23} \\ c_{31} & c_{32} & c_{33}\end{bmatrix} =  \begin{bmatrix} 1 & 1 & 1 \\ x_1 & x_2 & x_3 \\ y_1 & y_2 & y_3\end{bmatrix}^{-1}
\end{equation}
is also linear. This is shown in Fig. \ref{triangle}, where the direct mapping of a grid of points on the $Z$-domain (small black dots, left plot) are mapped to the $W$-domain (red cicles, right plot) and then mapped back to the $Z$-domain (red cicles, left plot). 

Note that, in the polynomial mapping the linearity holds in all the cases when the number of control points is equal to the dimensions of the mapping domain plus one. However, in the more interesting general case, when the number of control points is greater, then the mapping becomes nonlinear, with no closed-form inverse mapping expression. However, the mapping becomes more flexible to describe complex domains. This is shown in the next subsection where the domain is identified by 4 points.

\subsection{4-points polynomial mapping}

The polynomial mapping using $n = 4$ points on bi-variate domains (quadrangular polygons) is nonlinear. However, an inverse mapping is still possible to derive, as shown in subsection \ref{fpi}. This means that this non-linear mapping is bijective.

\begin{figure}[ht]
	\centering\includegraphics[width=0.8\linewidth]{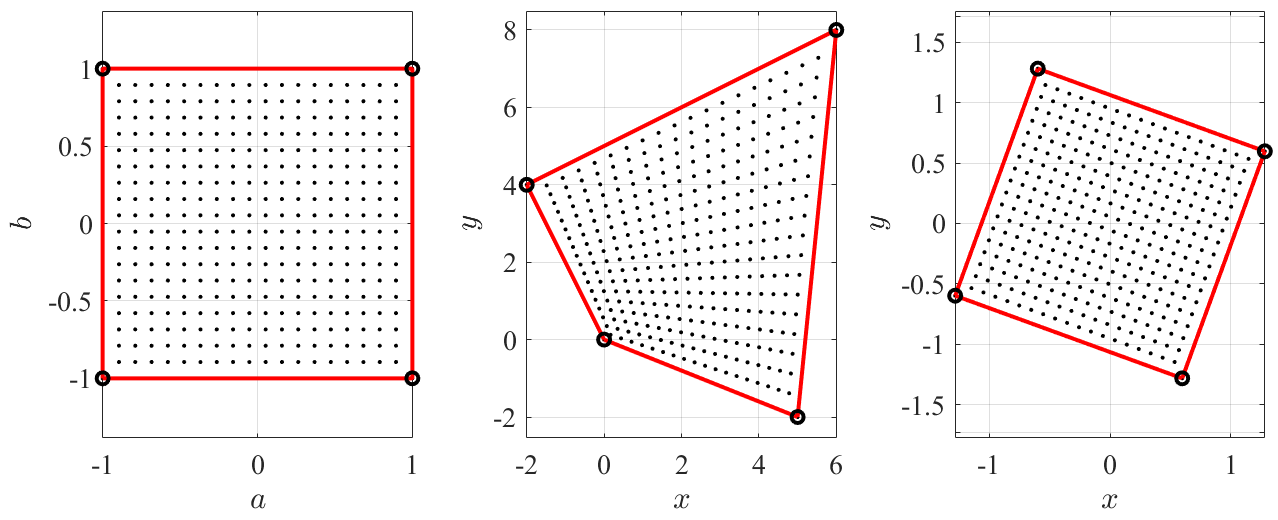}
	\caption{Polynomial mapping examples with 4 control points}
	\label{Fig1}
\end{figure}

Again, the mapping is performed between the unit-square domain, $(a, b) \in (-1, +1) \times (-1, +1)$, shown in top-left of Fig. \ref{Fig1}, and the quadrangular polygonal domain, defined by the four points, $x_k = [0, \, 5, \, 6, \, -2]$ and $y_k = [0, \, -2, \, 8, \, 4]$, shown in top-right of Fig. \ref{Fig1}. In particular, the mapping has been performed using a grid of $400$ points.

The direct mapping, $(a, b) \to (x,y)$, is,
\begin{equation}
x = \sum_{k = 1}^4 x_k \, f_k (a, b) \qquad \text{and} \qquad
y = \sum_{k = 1}^4 y_k \, f_k (a, b)
\end{equation}
where, the quadratic switching polynomial functions functions, $f_k, (a, b)$, have the following expressions,
\begin{equation}
\begin{cases}
f_1 (a, b) = \dfrac{1}{4} (1 - a - b + a \, b) \\[8pt]
f_2 (a, b) = \dfrac{1}{4} (1 + a - b - a \, b) \\[8pt]
f_3 (a, b) = \dfrac{1}{4} (1 + a + b + a \, b) \\[8pt]
f_4 (a, b) = \dfrac{1}{4} (1 - a + b - a \, b)
\end{cases}
\end{equation}

Figure \ref{Fig1} shows the direct mapping to the quadrangular polygonal domain defined by the four points (center figure), $x_k = [0, \, 1, \, 0, \, -1]$ and $y_k = [-1, \, 0, \, 1, \, 0]$, and to a domain obtained by rotating clockwise by $30$ deg. the $Z$ control points. This has been done using a grid of $400$ points.

Since this mapping transform lines to lines, the generic polygon side can be obtained simply by,
\begin{equation}
\B{w}_{ij} (t) = \B{w}_i + t \, (\B{w}_j - \B{w}_i) \qquad \text{where} \quad t \in [0, +1]
\end{equation}

\subsection{4-points inverse mapping}\label{fpi}

The inverse mapping, $(x,y) \to (a, b)$, can be obtained using the four boundary lines,
\begin{equation*}
\begin{cases}
\B{w}_{12} (t_1) = \B{w}_1 + t_1 \, (\B{w}_2 - \B{w}_1) \\
\B{w}_{43} (t_1) = \B{w}_4 + t_1 \, (\B{w}_3 - \B{w}_4)
\end{cases} \quad \text{and} \quad \begin{cases}
\B{w}_{14} (t_2) = \B{w}_1 + t_2 \, (\B{w}_4 - \B{w}_1) \\ 
\B{w}_{23} (t_2) = \B{w}_2 + t_2 \, (\B{w}_3 - \B{w}_2)
\end{cases}
\end{equation*}
to obtain the two lines passing through the point $(x, y)$,
\begin{equation*}
\begin{Bmatrix} x \\ y\end{Bmatrix} = \B{w}_{12} (t_1) + t_2 [\B{w}_{43} (t_1) - \B{w}_{12} (t_1)] = \B{w}_{14} (t_2) + t_1 [\B{w}_{23} (t_2) - \B{w}_{14} (t_2)]
\end{equation*}
or in a scalar form,
\begin{equation*}
t_1 \begin{Bmatrix} x_2 - x_1 \\ y_2 - y_1\end{Bmatrix} + 
t_2 \begin{Bmatrix} x_4 - x_1 \\ y_4 - y_1\end{Bmatrix} + t_1 t_2 \begin{Bmatrix} x_3 - x_4 + x_1 - x_2\\ y_3 - y_4 + y_1 - y_2\end{Bmatrix} = \begin{Bmatrix} x - x_1 \\ y - y_1\end{Bmatrix}
\end{equation*}
which can be written as,
\begin{equation*}
t_1 \begin{Bmatrix} C_{x1} \\ C_{y1}\end{Bmatrix} + t_2 \begin{Bmatrix} C_{x2} \\ C_{y2}\end{Bmatrix} + t_1 t_2 \begin{Bmatrix} C_{x3} \\ C_{y3}\end{Bmatrix} = \begin{Bmatrix} C_{x4} \\ C_{y4}\end{Bmatrix}
\end{equation*}
Setting,
\begin{equation*}
\begin{cases}
A = C_{x3} \, C_{y2} - C_{x2} \, C_{y3} \\ 
B = C_{x1} \, C_{y2} - C_{x2} \, C_{y1} + C_{x4} \, C_{y3} - C_{x3} \, C_{y4} \\
C =  C_{x4} \, C_{y1} - C_{x1} \, C_{y4}
\end{cases}
\end{equation*}
allows to write the solution in terms of a  quadratic equation,
\begin{equation*}
t_2 = \dfrac{-B + \sqrt{B^2 - 4 A C}}{2 \, A} \qquad \text{and} \qquad t_1 = \dfrac{x - x_1 -  C_{y2} \, t_2}{C_{y1} + C_{y3} \, t_2}
\end{equation*}
and finally,
\begin{equation*}
a = 2 \, t_1 - 1 \qquad \text{and} \qquad b = 2 \, t_2 - 1
\end{equation*}

\subsection{8-points polynomial mapping}

Constrained expressions for cubic quadrangular domains can be obtained by mapping the domain $(a, b) \in (-1, 1)\times(-1, 1)$ unit-square domain to a cubic quadrangular domain $(x, y)$, where each quadrangular side is identified by three points. This non-linear mapping is given by,
\begin{equation}
x = \sum_{k = 1}^8 x_k \, f_k (a, b) \qquad \text{and} \qquad
y = \sum_{k = 1}^8 y_k \, f_k (a, b)
\end{equation}
where, the cubic polynomials functions $f_k, (a, b)$ are the switching ($f_i (a_j, b_j) = \delta_{ij}$) functions,
\begin{equation}
\begin{cases}
f_1 (a, b) = \dfrac{1}{4} (a - 1) (1 - b) (a + b + 1) \\[8pt]
f_2 (a, b) = \dfrac{1}{2} (1 - a^2) (1 - b) \\[8pt]
f_3 (a, b) = \dfrac{1}{4} (a + 1) (1 - b) (a - b - 1) \\[8pt]
f_4 (a, b) = \dfrac{1}{2} (a + 1) (1 - b^2) \\[8pt]
f_5 (a, b) = \dfrac{1}{4} (a + 1) (b + 1) (a + b - 1) \\[8pt]
f_6 (a, b) = \dfrac{1}{2} (1 - a^2) (b + 1) \\[8pt]
f_7 (a, b) = \dfrac{1}{4} (a - 1) (b + 1) (a - b + 1) \\[8pt]
f_8 (a, b) = \dfrac{1}{2} (1 - a) (1 - b^2)
\end{cases}
\end{equation}
The quadrangular boundaries can be obtained by setting $a = \pm 1$ and $b = \pm 1$, respectively.
\begin{equation}\label{boundaryfunctions}
\begin{cases}
\B{c}_1 (a) = \dfrac{1}{2} (a - 1) a \, \begin{Bmatrix*}[r] -1 \\ -1\end{Bmatrix*} + (1 - a^2) \, \begin{Bmatrix*}[r] 0 \\ -1\end{Bmatrix*} + \dfrac{1}{2} (a + 1) a \, \begin{Bmatrix*}[r] 1 \\ -1\end{Bmatrix*} \\[7pt]
\B{c}_2 (b) = \dfrac{1}{2} (b - 1) b \, \begin{Bmatrix*}[r] 1 \\ -1\end{Bmatrix*} + (1 - b^2) \, \begin{Bmatrix*}[r] 1 \\ 0\end{Bmatrix*} + \dfrac{1}{2} (b + 1) b \, \begin{Bmatrix*}[r] 1 \\ 1\end{Bmatrix*} \\[7pt]
\B{c}_3 (a) = \dfrac{1}{2} (a - 1) a \, \begin{Bmatrix*}[r] -1 \\ 1\end{Bmatrix*} + (1 - a^2) \, \begin{Bmatrix*}[r] 0 \\ 1\end{Bmatrix*} + \dfrac{1}{2} (a + 1) a \, \begin{Bmatrix*}[r] 1 \\ 1\end{Bmatrix*} \\[7pt]
\B{c}_4 (b) = \dfrac{1}{2} (b - 1) b \, \begin{Bmatrix*}[r] -1 \\ -1\end{Bmatrix*} + (1 - b^2) \, \begin{Bmatrix*}[r] -1 \\ 0\end{Bmatrix*} + \dfrac{1}{2} (b + 1) b \, \begin{Bmatrix*}[r] -1 \\ 1\end{Bmatrix*}
\end{cases}
\end{equation}

\begin{figure}[!ht]
	\centering\includegraphics[width=\linewidth]{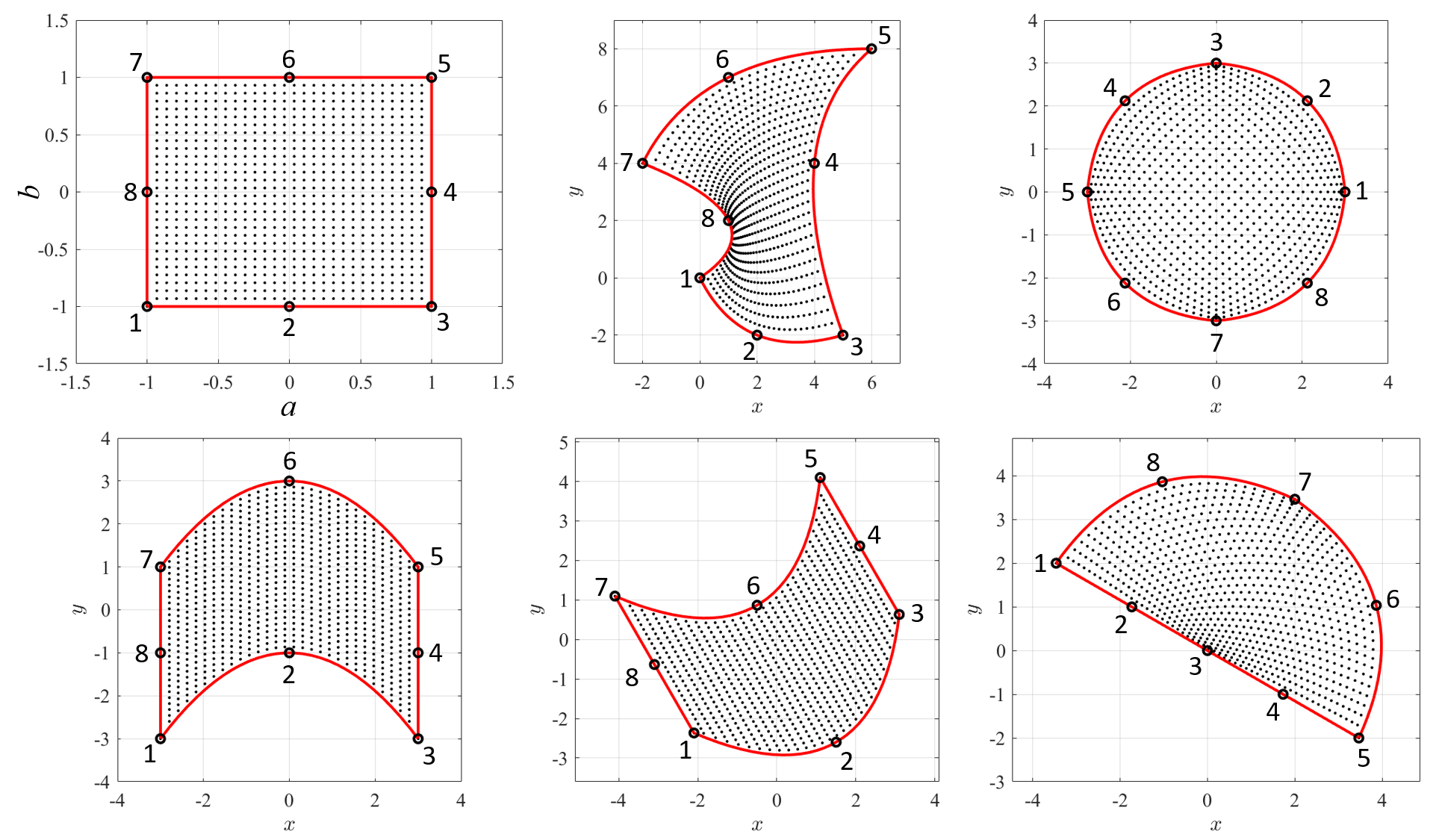}
	\caption{Polynomial mapping examples with 8 control points}
	\label{Fig2}
\end{figure}

Figure \ref{Fig2} shows the transformation of $900$ points from the $Z$ domain (top-left) to five different domains in $W$ defined by the boundary control points,
\begin{eqnarray}
\text{Top-center} & \to & \begin{cases}
\B{x}_k = \{0, 5, 6, -2, 2, 4, 2, 1\}\T \\
\B{y}_k = \{0, -2, 8, 4, -2, 4, 7, 2\}\T
\end{cases} \nonumber \\
\text{Top-right} & \to & \begin{cases}
\B{x}_k = 3 \, \cos((k - 1) \, \pi/4) \\
\B{y}_k = 3 \, \sin((k - 1) \, \pi/4)
\end{cases} \quad \text{where} \quad k = 1, \cdots, 8 \nonumber \\
\text{Bottom-left} & \to & \begin{cases}
\B{x}_k = \{-3, 0, 3, 3, 3, 0, -3, -3\}\T \\
\B{y}_k = \{-3, -1, -3, -1, 1, 3, 1, -1\}\T
\end{cases} \nonumber
\end{eqnarray}
\begin{eqnarray}
\text{Bottom-center} & \to & \text{Bottom-left points rotated by} \; 5\pi/6 \; \text{rad} \nonumber \\
\text{Bottom-right} & \to & \begin{cases}
\B{x}_k = \{-4, -2, 0, 2, 4, 4\cos(\pi/4), 0, -4\cos(\pi/4)\}\T \nonumber \\
\B{y}_k = \{0, 0, 0, 0, 0, 4\cos(\pi/4), 4, 4\cos(\pi/4)\}\T\end{cases} \nonumber
\end{eqnarray}
Specifically, the points in the top-right figure are selected as counter-clockwise points rotated by $45$ deg. The (cubic) boundary is not analytically circular, but it differ from it by roughly 1\%. Using more points a circular domain can be approximated with higher level of accuracy.

Unfortunately, the inverse of the polynomials mapping using $8$ points has not been found. Because of this, the following section provides a least-squares procedure to estimate an approximated inverse mapping for the 8 (or more) points polynomial mapping as well as for complex conformal mapping with no inverse mapping.

\section{Approximate least-squares inverse mapping}\label{ApproxInverse}

Inverse mappings always exist for the projection mapping, only. In fact and in general, complex and polynomial mappings do not admit closed-form inverse. However, when the direct mapping considered is bjiective, then an approximated least-squares inverse mapping is proposed in this section. The inverse mapping is needed because the existing TFC methodology has been developed for any dimensional space and for a wide set of constraint types, but for rectangular domains, only. Therefore, in order to apply the TFC framework to generic domains, an inverse mapping function is required to map the problem domain into the $Z$ rectangular domain where to apply the TFC framework.

To that end, a least-squares approximate inverse mapping is proposed as a linear combination of orthogonal polynomials,
\begin{equation}\label{InverseComplex}
\hat{z} = \xi_k \, \psi_k (w)
\end{equation}
where $\psi_k (w)$ is the $k$-th
orthogonal polynomial and where $z$ and $w$ are the complex coordinates that can be used, for real mapping, as $z = a + i \, b$ and $w = x + i\, y$. 

Let $\B{z}$ be the vector of $N$ discretized points (e.g., grid or uniformly distributed) in the $Z$ domain and $\B{w} = f (\B{z})$ be the direct mapping. Then, Eq. (\ref{InverseComplex}) can be specified for all the points and a least-squares solution of the unknown coefficients vector, $\B{\xi}$, can be estimated. It is important to outline that before applying this least-squares approach, translation and scaling might be necessary to contain all the $\B{w}$ points in a scaled domain in which the polynomial basis functions are defined as, for instance, in the $[-1, +1]$ range for Legendre or Chebyshev orthogonal polynomials.

\begin{figure}[!ht]
	\centering\includegraphics[width=\linewidth]{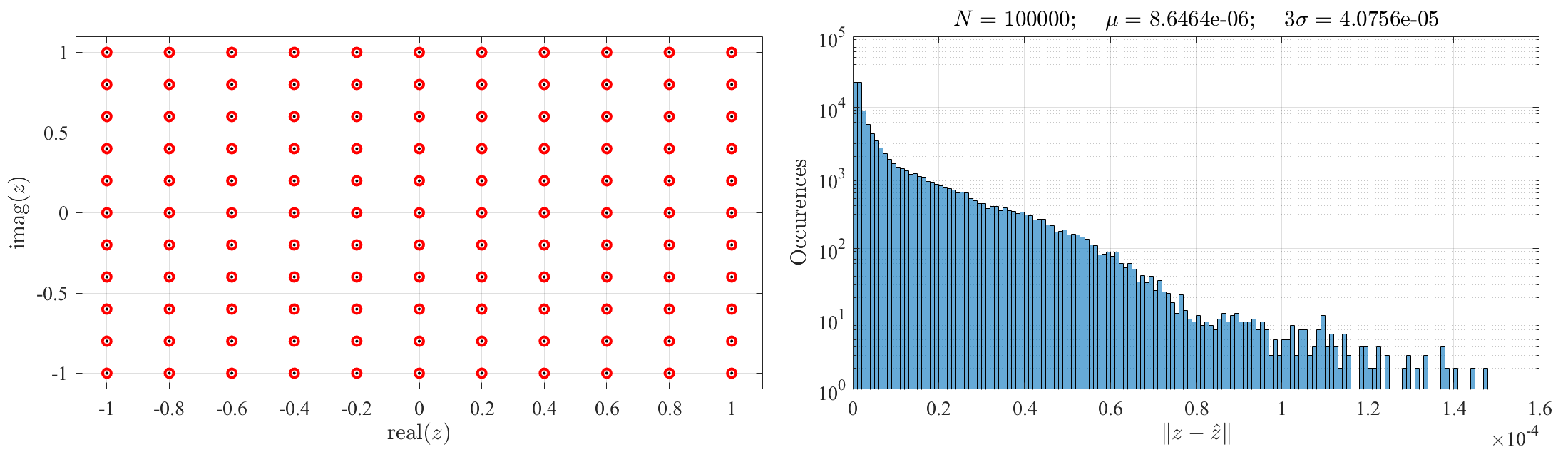}
	\caption{Approximate inverse function tests result}
	\label{Fig4}
\end{figure}

Figure \ref{Fig4} shows the accuracy results for the proposed approximate inverse mapping least-squares estimate. In this example, the same control points defined in Eq. (\ref{points}) have been used with $N = 121$ grid points in the $(-1, +1)\times (-1, +1)$ $Z$-domain. The left figure shows the $z_i$ grid of points as little black dots. The mapping points set, $[z_i, w_i]$, are then used to generate the approximate inverse mapping using Eq. (\ref{InverseComplex}). The inverse mapping estimates, $\hat{z}_i$, were obtained using $20$ Chebyshev orthogonal polynomial functions, $\psi_k (w)$. The $\hat{z}_i$ estimates are shown, as small red circles, in the left figure. 

To quantify the accuracy of the proposed method, an $100,000$ Montecarlo tests were made with (uniformly distributed) random points in the $Z$ domain. The histogram of the error, $|z_i - \hat{z}_i|$ is shown in the right figure. The mean error, $\mu = 8.6464 \times 10^{-6}$, and the $3\sigma = 4.0756 \times 10^{-5}$ error can be considered pretty good for most of the accuracy requirements of practical applications. Note that, the approximate inverse transformation must be computed just once.

\begin{figure}[!ht]
	\centering\includegraphics[width=0.75\linewidth]{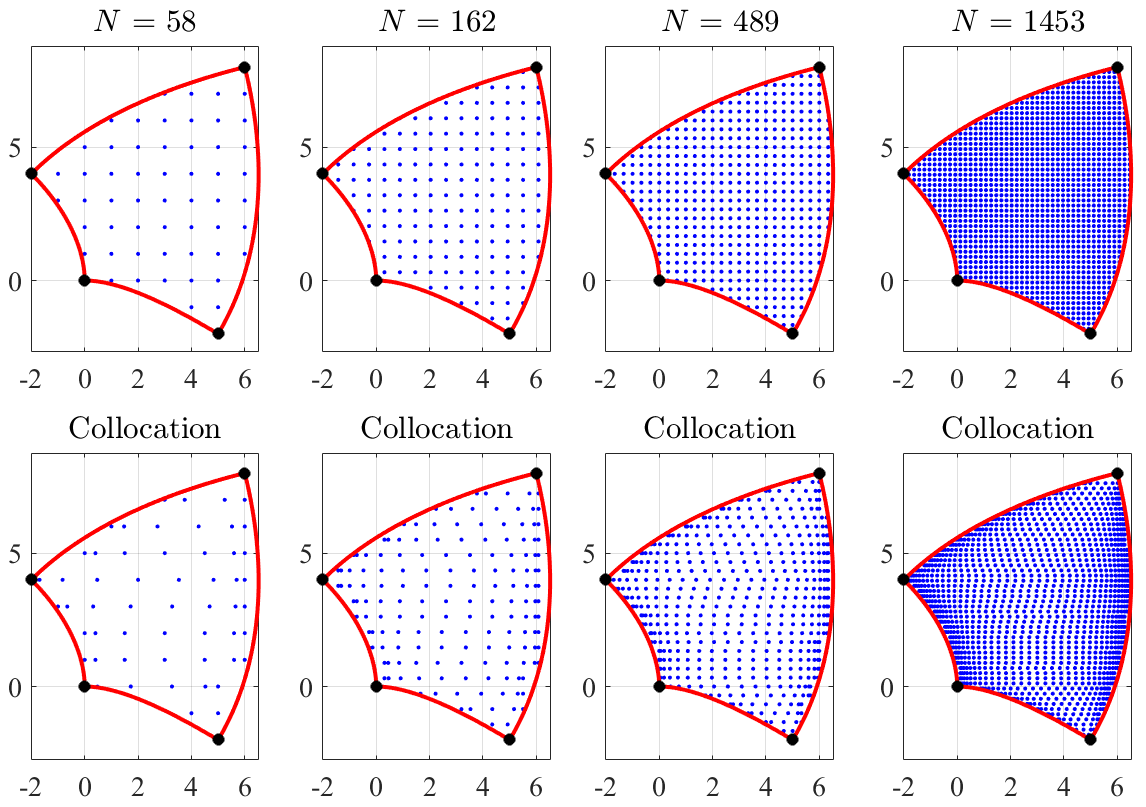}
	\caption{Points selection examples in the $W$ domain: grid (top) and collocation (bottom)}
	\label{grid}
\end{figure} 

Also, various distributions of points can be selected by embedding the $W$ domain into a slighlty larger rectangular domain. A grid distribution of points can be easily obtained for any convex $W$ domains and, using a more complex algorithm, for non-convex domains. 

Figure \ref{grid} shows, for $57$, $162$, $488$, and $1,453$ points, respectively, the selection of grid of points for the domain shown in Fig. \ref{Fig3}. Using a simple algorithm it is possible to obtain collocation-type of points. This is show in the bottom plots of Fig. \ref{grid} using the same number of points.

\section{Procedure to apply the \TFC\ to non rectangular domain}

This section summarizes, step-by-step, how to derive a constrained expression in the $W$ domain using the direct mapping function, $w = f (z)$. In this section $z$ and $w$ indicate either, the complex or the ($z \equiv [a,b]$ and $w \equiv [x,y]$) real coordinates in the $Z$ and $W$ domains.
\begin{enumerate}
	\item The $W$ domain (or subdomain, as those shown in Fig. \ref{multi}) is identified by four boundaries, $c_i (w) = 0$, along which the Dirichlet constraints are defined by the four functions, $v_i (w)$.
	\item A set of $N_g$ grid (or collocation-type) points within the $W$ domain, $w_k$, is computed (see Fig. \ref{grid}),
	\item If the inverse function, $z = f^{-1} (w)$, is available (as, for instance, when using the projection mapping or an invertable complex mapping) then the $w_k$ points are mapped back to the $Z$ domain, $z_k = f^{-1} (w_k)$. In the case the inverse mapping function is not available, then a least-squares approximated inverse function, $z = \hat{f}^{-1} (w)$, can be derived as described in Section \ref{ApproxInverse}.
	\item A set of points, $w_b$, belonging to the boundary, $c_i (w_b) = 0$, are computed. At these points/coordinates the Dirichlet constraint has value $v_i (w_b)$.
	\item Using the inverse (or the approximated inverse) mapping function, $z_b = f^{-1} (w_b)$, the boundary coordinates in the $Z$ domain are computed. At these corresponding points the constraint values are $v_i (w_b)$.
	\item The boundary constraints in the $Z$ domain ares obtained by interpolating the points $[z_b, \, v_i (w_b)]$ {\it and} using a constrained expression to guarantee the continuity at the four corners, as explained in \cite{TFC01}. 
	\item A Coons patch \cite{Coons} is then used to obtaing the simplest interpolating surface satisfying the Dirichlet boundary constraints in the $Z$ domain. This means obtaining the surface, $S_c (z)$, providing the value at coordinate $z$.
	\item A grid of points of the $S_c (z)$ surface is then mapped to the $W$ domain, $w_c = f (z_c)$
	\item Finally, the TFC functional in the $W$ domain is then given by,
	\begin{equation}
	\boxed{ S (w, g(w)) = S_c \left(f^{-1} (w)\right) + g (w) \ds\prod_{i = 1}^4 c_i (w) }
	\end{equation}
\end{enumerate}

\section{Conclusions}

This study consists of an initial investigation on how to apply the \TFC, which has been developed for rectangular domains in any dimensional space, to generic domains in 2-dimensional spaces. This has been done by three distinct approaches: 1) complex (conformal) mapping, 2) projection mapping, and 3) polynomial mapping. Discussions and examples are provided to highlight the features of each one of these three mappings, such as, conformal property, invertible property, or flexibility to describe different domains.

For the cases of bijective mappings where no analytical inversion is known, such as for some complex mapping and for the polynomial mapping, a method to derive least-squares approximate inverse mapping is provided. In addition, this study also describes how to replace constraint boundaries defined by a sequence of different functions by a single equation. This is required to derive a Coons patch in the $Z$ domain. 

This manuscript represents the first generalization of the \TFC\ for non rectangular domains. As such, there is additional research to be performed to better clarify details and to extend the methodology even more and to study its potential applications. For instance, given a generic $W$ domain (typically, polygonal or circular) of a real application, it would be interesting to derive a general methodology to obtain the mapping that approximates the $W$ domain boundaries by least-squares or by some other optimization techniques. Another possible future research includes the study of optimization techniques to search for the free function, $g (w)$, that allows to describe the set of four Dirichlet boundary constraints.

In other words, this study has the purpose to initiate the extension of the \TFC\ to non-rectangular domains.

\subsection*{The following abbreviations are used in this manuscript:}
\noindent 
\begin{tabular}{@{}ll}
	ODE & Ordinary Differential Equation\\
	PDE & Partial Differential Equation\\
	TFC & Theory of Functional Connections
\end{tabular}

\bibliographystyle{plain}

\end{document}